\documentclass[reqno,english]{amsart}

\usepackage[foot]{amsaddr}
\usepackage{amsmath,amsfonts,amssymb,graphicx,amsthm,enumerate, url, bm, bbm, upgreek, mathtools}
\usepackage[noadjust]{cite}
\usepackage[letterpaper]{geometry}
\usepackage[colorinlistoftodos]{todonotes}
\usepackage[colorlinks=true]{hyperref}
\usepackage[capitalize]{cleveref}

\usepackage{dja_macros}

\usepackage{matlab-prettifier}
\newcommand*{\GD}{\mathrm{GD}}

\title{A note on the capacity of the Binary perceptron}
\author[D.\ J.\ Altschuler \and K.\ Tikhomirov]{
        Dylan J. Altschuler \and 
        Konstantin Tikhomirov 
        }

\address{Dylan J. Altschuler, Department of Mathematical Sciences, Carnegie Mellon University.}
\address{Konstantin Tikhomirov, Department of Mathematical Sciences, Carnegie Mellon University.}
\begin{document}

\maketitle

\begin{abstract}
    Determining the capacity $\alpha_c$ of the Binary Perceptron is a long-standing problem. Krauth and M\'ezard (1989) conjectured an explicit value of $\alpha_c$, approximately equal to $.833$, and a rigorous lower bound
    matching this prediction was recently established
    by Ding and Sun (2019).
    Regarding the upper bound, Kim and Roche (1998) and Talagrand (1999)
    independently showed that $\alpha_c < .996$, while
    Krauth and M\'ezard outlined an argument which can be used to show that
    $\alpha_c < .847$. The purpose of this expository note is to record a complete proof of the bound $\alpha_c < .847$. The proof is a conditional first moment method combined with known results on the spherical perceptron. 
\end{abstract}


\section{Introduction} 
Let $A \in \RR^{\alpha N \times N}$ be a matrix with independent and identically distributed standard Gaussian entries. Consider the problem of finding $\sigma \in \cb{-1,+1}^N$ with $A\sigma > 0$, entrywise. The largest value of $\alpha$ such that this problem is solvable  is the random quantity $\alpha_* := \alpha_*(A)$, called the \textit{capacity}. The long-standing binary perceptron problem asks for the value $\ac$ that the capacity concentrates around. 

For intuition, observe $\ac \le 1$ is a trivial bound. Indeed, define the nested random sets $Z_t$,
\[
    Z_t := \cb{\sigma \in \cb{-1,+1}^N ~:~ A_i \sigma > 0, \quad \forall i \in 1,\dots,t}\,.
\]
For any fixed $\sigma \in Q$, we have $(A\sigma)_i > 0$ with probability one half, independently for each $i \in [t]$. So $\E{|Z_{t}|} = 2^{N - t}$, and thus for any $\eps$ by Markov's inequality: \[
    \PP{\alpha_* > 1+\eps} = \PP{| Z_{(1+\eps)N}| \ge 1}< 2^{-\eps N}\,.
\]
A long-standing and notoriously difficult open problem is to rigorously verify the following conjecture. 

\begin{conjecture}[Krauth and Mezard \cite{k-m}, 1989]
    The capacity of the binary perceptron concentrates around an explicit constant $\ac \approx .833$. 
\end{conjecture}
The value of $\ac$ is given as the fixed point of the so-called TAP, or Thouless-Anderson-Palmer, equations. This characterization is not yet known to be well-posed due to uniqueness issues. The TAP equations can be derived from non-rigorous but powerful heuristics such as the replica method or cavity method. Current progress can be summarized as follows: 

\begin{theorem}[Ding and Sun \cite{ding-sun}; Kim and Roche \cite{k-r}; Talagrand \cite{tala-rsa}; Xu \cite{xu}]
    With high probability, the capacity of the binary perceptron is at least $\alpha_c$ (i.e. $.833...$) and at most $.9963$.
\end{theorem}
 The lower bound was proven in a breakthrough of Ding and Sun. Their result is conditional on some numerical assumptions---this is not atypical in constraint satisfaction literature---and was originally proven to hold only with positive probability. The mismatch between positive probability and high probability was overcome by Xu \cite{xu} who showed the existence of a sharp threshold sequence (see also \cite{naka-sun} for a later, different proof, and \cite{dja-zero} for an alternative perspective). The upper-bound $.9963$ was shown by Kim and Roche, and a different proof (establishing an upper-bound of $1-\eps$ for some small, non-explicit $\eps > 0$) was subsequently discovered by Talagrand.
 In \cite[section 5]{k-m}, Krauth and Mezard outlined an argument
 which can be used to prove the following theorem:
\begin{theorem}\label{th: main}
    With high probability, the capacity of the binary perceptron is at most $.847$.
\end{theorem}
The proof of the theorem is based on a deep result of Shcherbina and Tirozzi \cite{sch-tir} about the ``spherical perceptron,'' a related model that is far better understood.
The paper \cite{k-m} appeared long before
the spherical perceptron was solved in \cite{sch-tir},
and for that reason the argument in \cite[section 5]{k-m} could not be made fully rigorous at the time. The purpose of our note is to provide a complete proof of Theorem~\ref{th: main} as
a corollary of the main result of \cite{sch-tir}.

\subsection{Spherical Perceptron} The closely related \textit{spherical} perceptron problem asks for the largest $\alpha$ so that there exists $\sigma \in S^{N-1}$ with $A\sigma > 0$ entrywise. The spherical perceptron predates the binary perceptron and was originally studied as a model of the fundamental limits of machine learning. It is also significantly better understood than the binary perceptron: since $\sigma$ is drawn from a convex set in the spherical case, one has access to a variety of powerful methods from convex geometry. 

Recall for the binary perceptron, one can track the size of $|Z_t|$, the number of remaining solutions after adding $t$ constraints. For the spherical perceptron it makes sense to track the volume of the solution set instead. Since volume scales exponentially in dimension, it is natural to actually consider the logarithm of the volume. Since this can possibly be negative infinity, it is further sensible to truncate the volume at some super-exponential cutoff. This motivates the study of the following random quantity, sometimes called the free energy:
\begin{equation}\label{eq:F}
    F(A) := \frac{1}{N}\max\cb{ \ln\pa{ \frac{\mathrm{Vol}\pa{ B(0,\sqrt{N}; \RR^N) \cap \cb{\sigma:A\sigma > 0}}}{\mathrm{Vol}\pa{ B(0,\sqrt{N}; \RR^N)}}}, -N^2  } \,.
\end{equation}
The celebrated prediction of Gardner \cite{gard1, gard2} (later in collaboration with Derrida  \cite{gard-derr}) from 1987, rigorously established by Shcherbina and Tirozzi in 2002 \cite{sch-tir}, gives an exact formula for the free energy $F$. 
\begin{theorem}[Gardner-Derrida formula for the Spherical Perceptron \cite{sch-tir}]\label{thm:sphere}
Let $F := F(A)$, where $A$ is an $\alpha N \times N$ matrix with iid standard Gaussian entries. 
\begin{itemize}
    \item $\Var{F} \to 0$ as $N \to \infty$.
    \item For any $\alpha > 2$, $F \to -\infty$ as $N \to \infty$. 
    \item For any $\alpha < 2$, letting $u \in \RR$ be distributed as a standard normal, we have:
    \begin{equation}
        \lim_{N \to \infty} \E{F} = \GD(\alpha) := \min_{q \in [0,1]} \alpha\, \E{\ln H\pa{\frac{u\sqrt{q}}{\sqrt{1-q}}}} + \frac{1}{2}\frac{q}{1-q} + \frac{1}{2}\ln(1-q)\,,
    \end{equation}
    where $H(\sigma) = \PP{Z \ge \sigma}$ for $Z$ a standard normal random variable. 
\end{itemize}    
\end{theorem}
We refer the reader to the monographs of Talagrand \cite{tal-mf1, tal-mf2} for detailed exposition. 

\subsection{Proof of Main Result} 
    Let $A$ be an $\alpha N \times N$ matrix, for $\alpha = .847$. We wish to show that with high probability, there does not exist $\sigma \in \cb{-1,+1}^N$ with $A\sigma > 0$ entrywise. As before, we proceed by Markov's inequality. We start with the following decoupling trick. Let $O$ be a random rotation uniform from the Haar measure. By basic symmetry considerations,
    \begin{equation}\label{eq:cond-expec}
        \E{|Z_{\alpha N}|} = \sum_{\sigma \in \cb{-1,+1}^N} \PP{A\sigma > 0} = \sum_{\sigma \in \cb{-1,+1}^N} \PP{A(O\sigma) > 0}\,.
    \end{equation}
    Observe that $O\sigma$ is uniformly distributed on the sphere of radius $\sqrt{N}$. Recalling \cref{eq:F} 
    \[
        \PP{A(O\sigma) > 0~|~ A} = \frac{\mathrm{Vol}\pa{ B(0,\sqrt{N}; \RR^N) \cap \cb{\sigma:A\sigma > 0}}}{\mathrm{Vol}\pa{ B(0,\sqrt{N}; \RR^N)}} \le e^{NF(A)} \,.
    \]    
    By \cref{thm:sphere}, the event $\mathcal{G} := \cb{A: |F(A) - \GD(\alpha)| < \eps}$ has probability tending to one for any $\eps > 0$. We claim that for $\eps$ sufficiently small (e.g. $\eps = .001$), it holds $\mathcal{G} \subset \cb{A: F(A) < - \log 2 - .001}$.
    \begin{proposition}\label{prop:GD}
        $\GD(.847) < -\log 2 - .002$. 
    \end{proposition}

    Assuming the proposition, the theorem readily follows. By Markov's inequality:
    \begin{align*}
        \PP{|Z_{\alpha N}| \ge 1} &\le \co{1} + \PP{|Z_{\alpha N}| \ge 1,~\mathcal{G}} \le \co{1} + \E{|Z_{\alpha N}| \ind(\mathcal{G}) }\,.
    \end{align*}
    But, the event $\mathcal{G}$ implies $F(A) < -\log(2) - .001$. Then \cref{eq:cond-expec} yields:
    \begin{align*}
        \E{|Z_{\alpha N}| \ind(\mathcal{G}) } &\le \sum_{\sigma \in \cb{-1,+1}^N} e^{-N \log 2 - .001N}\,.
    \end{align*}
    In total, $\PP{|Z_{\alpha N}| \ge 1}  = \co{1}$ where $\alpha = .847$. Equivalently stated, the capacity is at most $\alpha$ with high probability. All that remains is the proof of the proposition.

\begin{proof}[Proof of Proposition]
        Since $\GD(\alpha)$ is a given by a minimization problem and we seek to establish an upper-bound, it suffices to check any particular value of $q$. Define
    \[
        \GD(\alpha, q) := \alpha\, \E{\ln H\pa{\frac{u\sqrt{q}}{\sqrt{1-q}}}} + \frac{1}{2}\frac{q}{1-q} + \frac{1}{2}\ln(1-q) \,.
    \]
    Take $q = 1/2$. This choice is informed by numerically computing the true minimizer to be $q \approx .504$, and also by the fact that taking $q = 1/2$ trivializes the computation of $\GD$. Indeed, letting $\phi$ and $\Phi$ denote the Gaussian PDF and CDF respectively, integration by parts yields:
    \[
        \E{\ln H(u)} = \int_{\RR} \phi(u) \log H(u) du= -\int_\RR (1-\Phi(u)) \frac{\phi(u)}{1-\Phi(u)} du = -1\,.
    \]
    Thus,
    \begin{align*}
        \GD(.847, .5) &:= .847~ \E{\ln H\pa{u}} + \frac{1}{2}\pa{1 + \ln\pa{\frac{1}{2}}} \\
        &= -.847 + \pa{.15342\dots} \\
        &< -\ln(2) - .002\,.
    \end{align*}
    \end{proof}    

\bibliography{perceptron}
\bibliographystyle{acm}
\newpage
\section*{Appendix}
The following is MATLAB code to numerically optimize the Gardner-Derrida formula. It is not explicitly used in the main body of this article. Nonetheless, we offer it to the interested reader since it was helpful for exploration. \\

\begin{lstlisting}[style=Matlab-editor]
    % Grid search over possible values of the order parameter "q" and the constraint density "alpha".
    qlist = .001:.001:.999;
    a = .846:.00005:.847;
    
    % Each possible value of "q" certifies an upper-bound on the Gardner-Derrida free energy formula. Track all "candidate" upper bounds.
    GDcand = zeros(length(alpha),length(qlist));

    gpdf = @(u) exp(-u.^2./2) ./ sqrt(2*pi)  ; % Gaussian PDF
    H = @(x) (1 - erf(x ./sqrt(2))) ./ 2;   % Complement of Gaussian CDF

    % evaluate GD(alpha, q). 
    for i = 1:length(qlist)
        q = qlist(i);

        % in order to avoid log(0), add a perturbation "eps" (of size 2^{-52}) to the argument for numerical stability.
        I = @(u) gpdf(u) .* log(eps('double') +   H(u .* sqrt(q) ./ sqrt(1-q))   );

        % array of candidate solutions (changed logarithm to base 2)
        GDcand(:,i) = (alpha .* integral(I,-10,10) + q./(2.*(1-q)) + log(1-q) ./ 2) ./ log(2);
    end
    close all
    GD = min(GDcand,[],2) ; % returns col vector GD(alpha)
    plot(alpha, GD)
    
    xlabel('\alpha')
    ylabel('GD')
    yline(-1, 'r', 'GD = -1') 
    xline(.84655, 'k', '\alpha = .84655')
\end{lstlisting}

\end{document}